\RequirePackage{fix-cm}
\documentclass[smallextended]{svjour3}       
\smartqed  
\usepackage{graphicx}
\usepackage{amsfonts}
\usepackage{mathrsfs}
\usepackage{amsmath}
\usepackage [latin1]{inputenc}
\begin{document}

\title{Approximating continuous function on orbit spaces
}


\author{Qianqian Xia          
}


\institute{Qianqian Xia \at
            Nanjing University of Information Science \& Technology, Nanjing, 210044, China\\
              \email{inertialtec@sina.com}
}

\date{Received: date / Accepted: date}

\maketitle

\begin{abstract}
In this paper we study a subclass of subcartesian space-the orbit space of a proper action of Lie group on smooth manifold. We show that continuous functions on orbit space can be approximated by smooth functions.
\end{abstract}
\section{Introduction}
It is well-known that a continuous function on a smooth manifold can be approximated by a smooth function $^{[1]}$ as stated by the following theorem:

{\bf Theorem 1.1}
 Let $M$ be a smooth manifold and $f: M \rightarrow \mathbb{R}$ be a continuous function. Then for any $\delta >0$, there exists smooth function $g\in C^{\infty}(M, \mathbb{R})$, such that
\[|g(x)-f(x)|< \delta,\]
for all $x\in M$.

There has long been perceived the need for an extension of the framework of
smooth manifolds $^{[2]}$  in differential geometry, which is too restrictive and does not admit
certain basic geometric intuitions. Sikorski¡¯s theory of differential spaces studies
the differential geometry of a large class of singular spaces which both contains the theory
of manifolds and allows the investigation of singularities. It is the investigation of
geometry in terms of differentiable functions: a differential structure on a topological space $S$ is a family
$C^{\infty}(S)$ of real-valued functions on $S$ satisfying three conditions.  Functions $f\in C^{\infty}(S)$ are called smooth functions on $S$.
A differential space $S$ is said to be subcartesian $^{[3]}$ if every point of $S$ has a neighbourhood
diffeomorphic to a subset of some Cartesian space $\mathbb{R}^n$.
The theory of subcartesian spaces has been developed by \'Sniatycki et.al.
in recent few years.  See ${[3]}$ for a systematic treatment on this topic.

Generally speaking, it is not valid that  a continuous function on a subcartesian space $S$ can be approximated by smooth functions on $S$. However,
we study a subclass of subcartesian space in this paper. It has been shown in ${[3]}$ that the
orbit space of a proper action of a connected Lie group on a smooth manifold is a subcartesian
space. In this paper we show that a continuous function on the orbit space $R$ can be approximated by smooth functions on it.

This paper is organized as follows. In Section 1, we recall some basic definitions on subcartesian space.
In Section 2, we recall some basic facts on orbit space. In Section 3, we prove our main results.

\section{Subcartesian space}
{\bf Definition 2.1}
A differential structure on a topological space $S$ is a family $C^{\infty}(S)$ of real-valued functions on $S$ satisfying the following conditions:

1. The family
\[\{f^{-1}(I)|f\in C^{\infty}(S) \text{and} \, I \, \text{is an open interval in} \, \mathbb{R}\}\]
is a subbasis for the topology of $S$.

2. If $f_1, \cdots, f_n \in C^{\infty}(S)$ and $F \in C^{\infty}(\mathbb{R}^n)$, then $F(f_1, \cdots, f_n) \in C^{\infty}(S)$.

3. If $f: S\rightarrow \mathbb{R}$ is a function such that, for every $x \in S$, there exist an open neighborhood $U$ of $x$,
and a function $f_x \in C^{\infty}(S)$ satisfying
\[f_x|U=f|U,\]
then $f \in C^{\infty}(S)$. Here, the subscript vertical bar $|$ denotes a restriction.
\quad

{\bf Definition 2.2}
A map $\phi: R \rightarrow S$ is $C^{\infty}$ if $\phi^*f=f\circ \phi \in C^{\infty}(R)$ for
every $f \in C^{\infty}(S)$. A $C^{\infty}$ map $\phi$ between differential spaces is a $C^{\infty}$ diffeomorphism
if it is invertible and its inverse is $C^{\infty}$.
\quad

An alternative way of constructing a differential structure on a set $S$, goes as follows. Let $\mathcal{F}$ be a family of real-valued
functions on $S$. Endow $S$ with the topology generated by a subbasis
\[\{f^{-1}(I)|f\in \mathcal{F}\, \text{and} \,I \,\text{is an open interval in}\, \mathbb{R}\}.\]
Define $C^{\infty}(S)$ by the requirement that $h \in C^{\infty}(S)$ if, for each $x\in S$, there
exist an open subset $U$ of $S$, functions $f_1, \cdots, f_n\in F$, and $F \in C^{\infty}(\mathbb{R}^n)$ such
that
\[h|_U = F(f_1,\cdots, f_n)|_U.\]
Clearly, $\mathcal{F}\in C^{\infty}(S)$. It is proved in [3] that $C^{\infty}(S)$ defined
here is a differential structure on $S$. We refer to it as the differential structure
on $S$ generated by $\mathcal{F}$.

Let $R$ be a differential space with a differential structure $C^{\infty}(R)$, and let $S$
be an arbitrary subset of $R$ endowed with the subspace topology (open sets in
$S$ are of the form $S \cap U$, where $U$ is an open subset of $R$). Let
\[R(S)=\{f|_S|f \in C^{\infty}(R)\}.\]

{\bf Definition 2.3}
The space $R(S)$ of restrictions to $S \subseteq R$ of smooth
functions on $R$ generates a differential structure $C^{\infty}(S)$
on $S$ such that the
differential-space topology of $S$ coincides with its subspace topology. In this
differential structure, the inclusion map $i: S\rightarrow R$ is smooth.
\quad

In other words, $R(S)$ is the space of restrictions to $S$ of smooth functions on $R$.

Now consider an equivalence relation $\sim$ on a differential space $S$ with differential
structure $C^{\infty}(S)$. Let $R=S/{\sim}$ be the set of equivalence classes
of $\sim$, and let $\rho : S\rightarrow R$ be the map assigning to each $x\in S$ its equivalence
class $\rho(s)$.

{\bf Definition 2.4}
The space of functions on $R$, given by
\[C^{\infty}(R)=\{f: R\rightarrow \mathbb{R}|\rho^*f\in C^{\infty}(S)\},\]
is a differential structure on $R$. In this differential structure, the projection map
$\rho: S\rightarrow  R$ is smooth.
\quad

It should be emphasized that, in general, the quotient topology of  $R=S/{\sim}$
is finer than the differential-space topology defined by $C^{\infty}(R)$.

{\bf Definition 2.5}
A differential space $S$ is subcartesian if it is Hausdorff and every point $x \in S$ has a  neighbourhood $U$ diffeomorphic to a subset $V$ of $\mathbb{R}^n$ by the map $\phi: U\rightarrow \mathbb{R}^n$. $(U, \mathbb{R}^n, \phi)$ is said to be a local chart of $S$.  The subcartesian space
is said to be locally connected if for every point $x \in S$, there exists a connected neighborhood $U$ such that $(U, \mathbb{R}^n, \phi)$ is a local chart of $S$.
\quad

\section{Orbit space}
Consider a smooth and proper action
\begin{eqnarray}
\Phi:&&G\times M \rightarrow M\nonumber\\
&&(g,x)\rightarrow \Phi(g,x)=\Phi_g(x)=gx
\end{eqnarray}
of a Lie group $G$ on a manifold $M$.

We endow the orbit space $R=M/G$ with the quotient topology.
In other words, a subset $V$ of $R$ is open if $U=\rho^{-1}(V)$ is open in $M$, where $\rho:M\rightarrow R$ is the
canonical projection (the orbit map). Let
\[C^{\infty}(R)=\{f: R\rightarrow \mathbb{R}|\rho^*f\in C^{\infty}(M)\}.\]
$C^{\infty}(R)$ is a differential structure on $R$.

{\bf Proposition 3.1}$^{[3]}$
 The topology of $R$ induced by $C^{\infty}(R)$ coincides with the
quotient topology.

{\bf Proposition 3.2}$^{[3]}$
A slice through $x\in M$ for an action of $G$ on $M$ is a submanifold $S_x$ of $M$ containing $x$ such that:\\
1. $S_x$ is transverse and complementary to the orbit $Gx$ of $G$ through $x$. In other words,
\[T_xM=T_xS_x\oplus T_x(Gx).\]\\
2. For every $x'\in S_x$, the manifold $S_x$ is transverse to the orbit $Gx'$; that is,
\[T_{x'}M=T_{x'}S_x+T_{x'}(Gx').\]
3. $S_x$ is $G_x$-invariant.\\
4. Let $x'\in S_x$. If $gx' \in S_x$, then $g \in G_x$.
\quad

The existence of a slice through $x \in M$ is ensured by the following result.

{\bf Proposition 3.3}$^{[3]}$
There is an open ball $B$ in $\text{hor}T_xM$ centred at $0$ such that $S_x=\exp_x(B)$
is a slice through $x$ for the action of $G$ on $M$, where $\exp_x(v)$ is the value
at 1 of the geodesics of $G$-invariant Riemannian metric originating from $x$ in the direction $v$.
Further, the set $GS_x=\{gq | g\in G\, \text{and}\, q\in S_x\}$
is a $G$ invariant open neighbourhood of $x$ in $M$.
\quad

By construction, $S_x=\exp_x(B)$, where $\exp_x$ is an $H$-equivariant map from
a neighbourhood of $0$ in $T_xM$ to a neighbourhood of $x$ in $M$, and $B$ is a ball in
$\text{hor} T_xM$ invariant under a linear action of $H$ centred at the origin. The action of $H$ on $T_xM$ is linear, and it leaves
$\text{hor}T_xM$ invariant. Hence, it gives rise to a linear action of $H$ on $\text{hor}T_xM$.
Moreover, the restriction of $exp_x$ to $B$ gives a diffeomorphism $\psi : B\rightarrow S_x$,
which intertwines the linear action of $H$ on $\text{hor}T_xM$ and the action of $H$ on
$S_x$.

Since $B$ is an $H$-invariant open subset of $\text{hor}T_xM$ and the action of $H$ on
$\text{hor}T_xS_x$ is linear, by a theorem of G.W. Schwarz, smooth $H$-invariant functions
on $S_x$ are smooth functions of algebraic invariants of the action of $H$ on
$\text{hor} T_xM$. Let $\mathbb{R}[\text{hor} T_xM]^H$ denote the algebra of $H$-invariant polynomials on
$\text{hor} T_xM$. Hilbert¡¯s Theorem ensures that $\mathbb{R}[\text{hor} T_xM]^H$ is finitely generated.
Let $\sigma_1, \cdots, \sigma_n$ be a Hilbert basis for $\mathbb{R}[\text{hor} T_xM]^H$  consisting of homogeneous
polynomials. The corresponding Hilbert map
\begin{eqnarray}\label{eq:2}
\sigma: \text{hor} T_xM\rightarrow \mathbb{R}^n: v \rightarrow \sigma(v) =(\sigma_1(v),\cdots,\sigma_n(v))
\end{eqnarray}
induces a monomorphism $\tilde \sigma: (\text{hor} T_xM)/H\rightarrow \mathbb{R}^n : Hv\rightarrow \sigma(v)$, where $Hv$
is the orbit of $H$ through $v\in \text{hor} T_xM$ treated as a point in $(\text{hor} T_xM)H$. Let
$Q$ be the range of $\sigma$. By the Tarski¨CSeidenberg Theorem, $Q$ is a semi-algebraic
set in $\mathbb{R}^n$.
Let
\begin{eqnarray}\label{eq:3}
\phi : (\text{hor} T_xM)/H \rightarrow Q\subseteq \mathbb{R}^n
\end{eqnarray}
be the bijection induced by $\tilde \sigma$. $\phi$ is a diffeomorphism $^{[3]}$.

Since $B$ is an $H$-invariant open neighbourhood of $0$ in $\text{hor} T_xS_x$, it follows
that $B/H$ is open in $(\text{hor} T_xM)/H$. Hence, $B/H$ is in the domain of the diffeomorphism
$\phi: (\text{hor}T_xM)/H\rightarrow Q$, which induces a diffeomorphism of $B/H$
onto $\phi(B/H)\rightarrow Q\subseteq \mathbb{R}^n$. Thus, $B/H$ is diffeomorphic to a subset of $\mathbb{R}^n$.
But $B/H$ is diffeomorphic to $S_x/H$, and $S_x/H$ is diffeomorphic to $GS_x/G$.
Therefore, $GS_x/G$ is diffeomorphic to a subset of $\mathbb{R}^n$.

{\bf Theorem 3.4}$^{[3]}$
The orbit space $R=M/G$ of a proper action of $G$ on $M$ with
the differential structure $C^{\infty}(R)$ is subcartesian.

\section{Approximating continuous function on orbit spaces}
{\bf Lemam 4.1}$^{[1]}$
Let $U$ be an open subset of $\mathbb{R}^m$ and let $K\subseteq U$ be a compact subset of $\mathbb{R}^m$. Then for any continuous function
$f\in C^{0}(U)$ and any $\delta > 0$, there exists a smooth function $g\in C^{\infty}(U)$,  such that
\[|g(y)-f(y)|< \delta,\]
for any $y \in K$.\quad

{\bf Lemma 4.2}
For each $y_0 \in R$, there exist a local neighborhood $V$ of $y_0$ and a compact subset $K \subseteq V$ satisfying that for any continuous function
$f\in C^{0}(V)$ and any $\delta > 0$, there exists a smooth function $g\in C^{\infty}(V)$, where $(V, C^{\infty}(V))$ is a differential subspace of $R$, such that
\[|g(y)-f(y)|< \delta,\]
for any $y \in K$.
\quad

Proof. Let  $x_0\in M$ such that $\rho(x_0)=y_0$. Let  $H$ be the isotropy group of $x_0$ and  $S_{x_0}=\exp_{x_0}(B)$ be a slice through $x_0\in M$, where $\exp_{x_0}$ is an $H$-equivariant map from
a neighbourhood of $0$ in $T_{x_0}M$ to a neighbourhood of $x_0$ in $M$, and $B$ is a ball in
$\text{hor} T_{x_0}M$ invariant under a linear action of $H$ centred at the origin. Let $K\subseteq B$ be a $H$ invariant compact subset of $T_{x_0}M$.

Then for any continuous function $f$ on $\rho(GS_{x_0})$, it follows that $\exp^{*}(\rho^*f|_{S_{x_0}})$ is a continuous function on $B$.
From Lemma 4.1 we know that any for $\delta > 0$ there exists a  smooth function $h\in C^{\infty}(B)$,  such that
\[|h(x)-\exp^{*}(\rho^*f|_{S_{x_0}})(x)|< \delta,\]
for any $x \in K$.\quad

Now consider the smooth function $h\circ \exp^{-1}$ on $S_x$, which satisfies that $|h\circ \exp^{-1}(x)-\rho^*f|_{S_{x_0}}(x)|< \delta,$
for any $x \in \exp(K)$, where $\exp(K)$ is an $H$-invariant compact subset in $M$. Since $H$ is compact, we may average $h\circ \exp^{-1}$ over $H$, obtaining a
$H$-invariant function
\[\tilde h=\int_{H}\Phi^*_{g}(h\circ \exp^{-1})d\mu(g),\]
where $d\mu(g)$ is the Haar measure on $H$ normalized so that $vol H = 1$.

The set $GS_{x_0}$ is a $G$-invariant open neighbourhood of $x_0$ in $M$. We can define
a $G$-invariant function $\tilde {f}_1$ on $GS_{x_0}$ as follows. For each $x''\in GS_{x_0}$, there exists
$g\in G$ such that $x''=gx'$ for $x'\in S_x$, and we set
\[\tilde {f}_1(x'')=\tilde h(x').\]
Since $\tilde h$ is $H$-invariant, the function $\tilde{f}_1$ is well defined on $GS_{x_0}$ and is $G$-invariant. Hence we have a smooth function $f_1 \in C^{\infty}(\rho(S_{x_0}))$ such that $\rho^*f_1=\tilde {f}_1$. Besides, for each $y\in \rho(\exp(K))$, we have
\begin{eqnarray}
|f_1(y)-f(y)|&=&|\tilde h(x)-\rho^*f(x)|=|\int_{H}(\Phi^*_{g}(h\circ \exp^{-1}))(x)d\mu(g)-\int_{H}(\Phi^*_{g}(\rho^*f))(x)d\mu(g)|\nonumber\\
&=&|\int_{H}(\Phi^*_{g}(h\circ \exp^{-1})-\Phi^*_{g}(\rho^*f))(x)d\mu(g)|\nonumber\\
&\leq &\int_{H}|(\Phi^*_{g}(h\circ \exp^{-1})-\Phi^*_{g}(\rho^*f))(x)|d\mu(g)\nonumber\\
&< &\int_{H}\delta d\mu(g)\nonumber\\
&=&\delta.
\end{eqnarray}
Hence for $y_0 \in R$ and for $x\in \rho^{-1}(y)$,there exist a local neighborhood $\rho(GS_x)$ of $y_0$ and a compact subset $\rho(K) \subseteq \rho(GS_x)$ satisfying that for any continuous function
$f\in C^{0}(\rho(GS_x))$ and any $\delta > 0$, there exists a smooth function $f_1\in C^{\infty}(\rho(GS_x))$, where $(\rho(GS_x), C^{\infty}(\rho(GS_x)))$ is a differential subspace of $R$, such that
\[|f_1(y)-f(y)|<\delta,\]
for any $y \in \rho(K)$. Hence the result follows immediately.\quad

{\bf Lemma 4.3}$^{[2]}$ Let $U, V$ be two open subsets of the smooth manifold $M$ satisfying that $\bar U$ is compact and $\bar U \cap \bar V=\emptyset$. Then there exist smooth function $f \in C^{\infty}(M)$ such that
\begin{eqnarray}
&&f(x)=1, x\in \bar U, \nonumber\\
&&0< f(x)< 1, x\in {M/{\bar V \cup \bar U}}\nonumber\\
&&f(x)=0, x\in \bar V.
\end{eqnarray}

{\bf Lemma 4.4}  Let $x\in M$ and let $0 \in W\subseteq V\subseteq U\subseteq B$ be $H$-invariant open subsets of $\text{hor}T_xM$ such that $\bar W\subseteq V$ and $\bar V$ are compact, where  $B$ is satisfies that $\exp(B)=S_x$.
Let $G$ be an open subset of $R$ and $(\rho(\exp(U)), \psi)$ be the local coordinate for $R$ induced by the Hilbert map (2).
Let $f: R \rightarrow \mathbb{R}$ be a continuous map satisfying that $f|_G$ is smooth. Then for any $\delta > 0$, there exists a continuous map $g: R\rightarrow \mathbb{R}$, such that

(1) $g(y)=f(y)$,  for any $y\in R/\rho(\exp(V))$;

(2) $g|_{G\cup \rho(\exp(W))}$  is smooth;

(3) $|g(y)-f(y)|<\delta$, for all $y \in R$,

where $\sigma$ is given by (2).

Proof. It  follows from Lemma 4.3 that there exists smooth function $\eta: \text{hor}_xM \rightarrow \mathbb{R}$ such that
\begin{eqnarray}
&&\eta(x)=1, x\in \bar W, \nonumber\\
&&0< \eta(x) < 1, x\in {V/{\bar W}}\nonumber\\
&&\eta(x)=0, x\in \text{hor}T_xM/{V},
\end{eqnarray}
which yields a smooth function $\eta \circ \exp^{-1}$ on $S_x$. Since $V, W$ are $H$-invariant, then by averaging $\eta \circ \exp^{-1}$ over $G_x$ we get a  $G_x$ invariant smooth function $\tilde \eta$ on $S_x$ satisfying that
\begin{eqnarray}
&&\tilde{\eta}(x)=1, x\in \exp(\bar W), \nonumber\\
&&0< \tilde{\eta}(x) < 1, x\in \exp({V/{\bar W}})\nonumber\\
&&\tilde{\eta}(x)=0, x\in S_x/\exp(V),
\end{eqnarray}
which can be extended to a smooth $G$ invariant function on $M$. Hence we get a function $\bar \eta \in C^{\infty}(R)$ satisfying that
\begin{eqnarray}
&&\bar {\eta}(y)=1, y\in \rho(\exp(\bar W)), \nonumber\\
&&0< \bar {\eta}(y) < 1, y\in \rho(\exp({V/{\bar W}}))\nonumber\\
&&\bar{\eta}(y)=0, y\in R/\rho(\exp(V)).
\end{eqnarray}
It follows from Lemma 4.2 that the function $f|_U$ can be approximated by smooth functions on $U$. That is, for any $\delta> 0$,  there exists smooth function $g_0\in C^{\infty}(\rho(\exp(U)))$ such that $|g_0(y)-f(y)|< \delta$, for $y \in \rho(\exp(\bar V))$.

Since $f=(1-\bar \eta)f+\bar \eta f$, we define
\[g=(1-\bar \eta)f+\bar \eta g_0.\]
Then the result follows immediately. \quad

{\bf Lemma 4.5}$^{[1]}$
Let $M$ be a second countable locally compact Hausdorff topological space. Then there exist countable many open sets $G_1, G_2, \cdots, G_k, \cdots$
satisfying

(1) $cl(G_j)$ is compact, $j=1,2, \cdots$;

(2) $cl(G_j)\subseteq G_{j+1}, j=1,2, \cdots$;

(3) $\cup G_j=\cup cl(G_j)=M,$\\
where $cl(G_j)$ denotes the closure of $G_j,j=1,2,\cdots$.
\quad

{\bf Lemma 4.6}
There exist locally finite open covers  $(U_j)_{j\in \mathbb{Z}_{>0}}, (V_j)_{j \in \mathbb{Z}_{>0}}, (W_j)_{j \in \mathbb{Z}_{>0}}$ of $R$ such that $cl(U_j)\subseteq V_j, cl(V_j)\subseteq W_j$, and $cl(U_j), cl(V_j), cl(W_j)$ are compact, for
each $j>0$, where $(W_j,\mathbb{R}^{n_j}, \phi_j)$ is a local chart of $R$ induced by the Hilbert map (2).

Proof. From Lemma 4.5 we know that there exist countable open sets $G_1, \cdots, G_k,\cdots$ satisfying conditions (1), (2) and (3) in Lemma 4.5.
It follows that $cl(G_h)/{G_{h-1}}$ is compact,  $G_{h+1}/{cl (G_{h-2})}$ is open and $cl(G_h)/{G_{h-1}}\subseteq G_{h+1}/{cl (G_{h-2})}$.  On the other hand we know the  local charts induced by the Hilbert map (2) of $R$ form an open cover of $R$. Then for $y\in cl({G}_h)/{G_{h-1}}$, there exist a local chart $(\mathscr{V}, \phi)$ of  $y$ induced by the Hilbert map (2).  Consider the $H$ invariant open set $\sigma^{-1}((G_{h+1}/{cl({G}_{h-2})}\cap \mathscr{V})$ in $\text{hor}T_xM$, where $\rho(x)=y$. There exists an open ball $B_{\epsilon}$ such that
 $\bar B_{\epsilon}\subseteq \sigma^{-1}((G_{h+1}/{cl({G}_{h-2})}\cap \mathscr{V})$ centered at $0$. Let  $\mathscr{W}=\sigma(B_{\epsilon})$. Hence $\mathscr{W}$ is an open subset containing $y$ such that $cl(\mathscr{W}) \subseteq (G_{h+1}/{cl({G}_{h-2})})\cap \mathscr{V}$.

It follows that

(1) $y\in W\subseteq (G_{h+1}/{cl({G}_{h-2})})\cap \mathscr{V}$;

(2) $\phi(y)=0$ and $\phi(W)=\phi(\mathscr{W})\subseteq \phi(\mathscr{V})$; \\
Since $cl(W)=\sigma(\bar B_{\epsilon})$ and $\bar B_{\epsilon}$ is compact, it follows that
$cl(W)$ is compact.

Let $\mathscr{W}_1=\sigma(B_{\epsilon_1})$ where $0< \epsilon_1< \epsilon$. Then $\mathscr{W}_1$ is an open set containing $y$ such that $cl(\mathscr{W}_1)\subseteq W$.
Denote by $V=\mathscr{W}_1$.   And  let $\mathscr{W}_2=\sigma(B_{\epsilon_2})$ where $0<\epsilon_2\textless \epsilon_1$. Then $\mathscr{W}_2$ be an open set containing $y$ such that $cl(\mathscr{W}_2)\subseteq V$. Denote by $U=\mathscr{W}_2$.  Then we have $cl(U) \subseteq V$ and $cl(V)\subseteq W$.

Since  ${cl({G}_h})/{G_{h-1}}$ is compact, there exist finitely many points $y_{h, 1}, y_{h, 2}, \cdots$, $y_{h, k_h} \in cl({G}_h)/{G_{h-1}}$, such that the corresponding open sets $U_{h,1}, U_{h,2}, \cdots, U_{h,k_h}$  form an open cover of $cl(G_h)/G_{h-1}$.
We claim that the corresponding open sets
\[(U_{1,1}, U_{1,2}, \cdots, U_{1,k_1}; U_{2,1}, U_{2,2}, \cdots, U_{2, k_2}; \cdots),\]
\[(V_{1,1}, V_{1,2}, \cdots, V_{1,k_1}; V_{2,1}, V_{2,2}, \cdots, V_{2, k_2}; \cdots),\]
\[(W_{1,1}, W_{1,2}, \cdots, W_{1,k_1}; W_{2,1}, W_{2,2}, \cdots, W_{2, k_2}; \cdots),\]
satisfy the conditions in the lemma.   We only need to prove the local finiteness of $(W_{i,j})$. Given $y\in R$, assume that
$y \in G_r$, then it follows from the above construction that there exist finite many $W_{i,j}$ that intersect $G_r$. In fact,
\[W_{i,j}\cap G_r=\emptyset, i\geq r+2, 1\leq j \leq k_i.\]
This completes the proof of the lemma.
\quad

{\bf Lemma 4.7} Let $f: R\rightarrow \mathbb{R}$ be a continuous function on $R$. Then for any $\delta> 0$, there exists $g\in C^{\infty}(R)$, such that
\[|g(y)-f(y)|< \delta,\]
for any $y \in R$.

Proof. From Lemma 4.6 we know that there exists  locally finite open covers  $(U_j)_{j\in \mathbb{Z}_{>0}}$, $(V_j)_{j \in \mathbb{Z}_{>0}}$, $(W_j)_{j \in \mathbb{Z}_{>0}}$ of $R$ such that $cl(U_j)\subseteq V_j, cl(V_j)\subseteq W_j$, and $cl(U_j), cl(V_j), cl(W_j)$ are compact, for
each $j>0$, where $(W_j,\mathbb{R}^{n_j}, \phi_j)$ is a local chart of $R$ induced by the Hilbert map (2).

Set $W_0=\emptyset, f_0=f$. Assume we have continuous function $f_k$ on $R$ such that $f_k|_{G_k}$ is smooth, where
\[G_k=\cup_{j=0}^kW_j.\]
Then it follows from Lemma 4.4 that there exists continuous function $f_{k+1}$ on $R$, such that $f_{k+1}|_{G_{k+1}}$ is smooth, where
\[G_{k+1}=\cup_{j=0}^{k+1}W_j.\]
Besides, $f_{k+1}|_{M/V_{k+1}}=f_k$, and
\begin{eqnarray}
|f_{k+1}(y)-f_k(y)|<\frac{\delta}{2^{k+1}},
\end{eqnarray}
for all $y \in R$.

Let \[g(y)=\lim_{k\rightarrow \infty}f_k(y).\]
It follows from (9) that $g$ is well-defined. And
\[|g(y)-f(y)|< \delta,\]
for any $y \in R$.

We claim that $g \in C^{\infty}(R)$. For $y \in R$, there exists $l \in \mathbb{Z}_{> 0}$, such that  $y \in W_l$. Now consider the function $f_{l}, f_{l+1}, \cdots$ which are smooth on $W_j$. It follows that $\rho^*f_l, \rho^*f_{l+1}, \cdots$ are smooth functions on the open subsets $\rho^{-1}(W)$ of $M$, which satisfies
\begin{eqnarray}
|\rho^*f_{k+1}(x)-\rho^*f_k(x)|<\frac{\delta}{2^{k+1}},
\end{eqnarray}
for $k \geq l$ and  $x \in \rho^{-1}(W)$. It follows that $\lim_{k\rightarrow \infty}\rho^*f_k \in C^{\infty}(\rho^{-1}(W))$. Besides, it is also $G$ invariant which descends to $g$ on $W$. Hence $g|_W \in C^{\infty} (W)$. Since $y$ is arbitrary, we have $g\in C^{\infty}(R)$. This completes the proof.

{\bf Theorem 4.8} Let $f: R\rightarrow \mathbb{R}^n$ be a continuous function on $R$. Then for any $\epsilon > 0$, there exists $g\in C^{\infty}(R; \mathbb{R}^n)$, such that
\[|g(y)-f(y)|< \epsilon,\]
for any $y \in R$. Besides, $g$ is homotopic to $f$.

Proof. Let $\delta=\epsilon/{\sqrt n}$. It follows from Lemma 4.7 that there exist smooth functions $g_1, \cdots, g_n\in C^{\infty}(R)$ such that \[|g_i(y)-f_i(y)|< \delta,\]
for any $y \in R$, where $f=(f_1, \cdots, f_n)$. Consider the smooth map $g=(g_1, \cdots, g_n)$. We have that
\[|g(y)-f(y)|< \epsilon,\]
for any $y \in R$. Besides, define
\[F(t.y)=(1-t)g(y)+tf(y),\]
for $(t,y)\in I \times R$. It is obvious that $F$ defines a homotopy from $g$ to $f$. Hence the result follows immediately. \quad


\end{document}